\documentclass[pdflatex]{article}
\usepackage{amssymb,amsfonts,amsmath,latexsym,amscd}
\usepackage[all,cmtip]{xy}
\usepackage{url}
\usepackage{mathrsfs}
\usepackage{xcolor}
\usepackage{graphicx}

\usepackage{cases}
\usepackage{diagbox}
\usepackage{tikz}

\newtheorem{theorem}{Theorem}[section]
\newtheorem{lemma}[theorem]{Lemma}
\newtheorem{corollary}[theorem]{Corollary}
\newtheorem{proposition}[theorem]{Proposition}

\newtheorem{definition}[theorem]{Definition}

\newenvironment{proof}{{\par\addvspace{0.1cm}\noindent \bf Proof. }}{\hfill$\Box$\par\medskip}

\newtheorem{remark}[theorem]{Remark}

\numberwithin{equation}{section}

\def\a{\alpha}
\def\RR{\mathbb{R}}
\def\QQ{\mathbb{Q}}

\def\e{\varepsilon}
\def\de{\delta}
\def\la{\lambda}
\def\si{\sigma}
\def\ga{\gamma}
\def\Ga{\Gamma}

\def\vect#1{\mbox{\boldmath $#1$}}

\def\jb#1{\textcolor{black}{{#1}}} 
 
\def\j#1{\textcolor{black}{{#1}}}

\usetikzlibrary{arrows.meta, calc}

\begin{document}

\title{Distinguishing finite metric spaces via similarity spectra}

\author{Jun O'Hara\footnote{Supported by JSPS KAKENHI Grant Number 23K03083.}}
\maketitle

\begin{abstract}
We study spectra and characteristic polynomials of similarity matrices associated with finite metric spaces, where the similarity matrix of a finite metric space $X=\{x_1,\dots,x_n\}$ is given by 
$\displaystyle Z(q)=(q^{d(x_i,x_j)})_{i,j}$. 
We introduce two spectral invariants of finite metric spaces,
the $q$-spectrum and the transition $q$-spectrum,
defined respectively from $Z(q)$ and its transition matrix.
In the case of graphs, these invariants recover the adjacency spectrum and the Laplacian spectrum in the limit $q\to0$.

Our main result shows that the $q$-spectrum determines a large class of finite metric spaces under a natural nondegeneracy condition.
We also prove that all four-point metric spaces are determined by their $q$-spectra.
The key observation is that the coefficients of the characteristic polynomial of $Z(q)$ encode cycle structures of the underlying metric space. 
We further investigate the transition $q$-spectrum \jb{and show that strongly regular graphs with the same parameters have identical $q$-spectra and transition $q$-spectra, providing infinitely many non-isomorphic examples that cannot be distinguished by these invariants. Finally, we present computational examples comparing these invariants with classical graph spectra.}

\end{abstract}



\setcounter{tocdepth}{3}

%!!!!!!!!!!!!!!!!!!!!!!!!!!!!!!!!!!!!!%%%%%%%%%%%%%%%%%%%%%%%%%%%%%%%%%%%%
\section{Introduction}
%!!!!!!!!!!!!!!!!!!!!!!!!!!!!!!!!!!!!!%%%%%%%%%%%%%%%%%%%%%%%%%%%%%%%%%%%%
% 
Spectral invariants of matrices associated with combinatorial or geometric structures play a central role in linear algebra and graph theory.
Typical examples include the adjacency spectrum and the Laplacian spectrum of graphs.
In this paper, we study analogous spectral invariants associated with finite metric spaces.
Let
\[
X=\{x_1,\dots,x_n\}
\]
be a finite metric space.
We associate with $X$ the similarity matrix
\[
Z(q)=(q^{d(x_i,x_j)})_{i,j},
\]
where $q$ is a formal parameter and $d(x_i,x_j)$ denotes the distance between $x_i$ and $x_j$.
We study the characteristic polynomial and spectra of these matrices.

The first invariant introduced in this paper is the $q$-spectrum,
defined as the spectrum of $Z(q)$.
The second is the transition $q$-spectrum,
defined as te spectrum of the transition matrix associated with $Z(q)$ as the stochastic normalization that is obtained by dividing each row by its sum. 
These invariants extend classical graph spectra:
in the case of graphs,
the adjacency spectrum and Laplacian spectrum are recovered in the limit $q\to0$.
Therefore, the two invariants considered in this paper provide extensions of the two fundamental graph spectra, the adjacency spectrum and the Laplacian spectrum, to the setting of finite metric spaces.

Our main interest is the reconstruction problem:
to what extent can a finite metric space be recovered from these spectral invariants?
In the setting of finite metric spaces, this is closely related to the unassigned distance geometry problem (uDGP) and the homometric problem.
The key observation is that the coefficients of the characteristic polynomial of $Z(q)$ encode combinatorial information on lengths of cycles.
Using this interpretation,
we prove reconstruction theorems for broad classes of finite metric spaces.

Our main result shows that the $q$-spectrum determines a large class of finite metric spaces. To be precise, our theorem applies to metric spaces that satisfy the following non-degeneracy condition: among the $N_2={{n}\choose2}$ pairwise distances, the sums of any three distinct distances are all different. 
The set of metric spaces satisfying this condtion is open dense in the moduli space $\mathcal{M}_n$ by Roff-Yoshinaga (\cite{RY}), which is obtained by parametrizing pairwise distances subject to the triangle inequalities and taking the quotient by the natural action of the symmetric group, equipped with the quotient topology which is induced from the standard Euclidean topology.
The condition of our theorem implies in particular that all pairwise distances are distinct, and hence the symmetric group acts freely. As a result, the subset of $\mathcal{M}_n$ consisting of metric spaces satisfying our non-degeneracy condition avoids the singular locus of the quotient and inherits the local topology of the ambient Euclidean space. 

The key is that the characteristic polynomial coefficients encode cycle-length data, in particular, the pairwise distances and triangle perimeters. 
Moreover, we also show that all metric spaces on at most four points are completely determined by this invariant.

The reason why such a condition is required is as follows.
When constructing invariants determined by distances from a given space, the presence of intermediate symmetries may cause the data to become entangled, making it impossible to recover the original structure from the invariant.
This can be avoided when the symmetry is minimal, that is, when the space is sufficiently generic.
This phenomenon is also observed when attempting to reconstruct spaces using other invariants.

Our second result shows that the transition $q$-spectrum determines a finite metric space if the multiset of pairwise distances is linearly independent over $\QQ$. 
A finite metric satisfying this condition is generic in the sense that the complement of the set of such spaces is measure zero in $\mathcal{M}_n$. 
Interestingly, while the theoretical conditions for the transition $q$-spectrum are strictly stronger, our computational experiments reveal that it actually exhibits significantly higher distinguishing power in practice. 
This phenomenon is consistent with known \jb{tendency} in spectral graph theory, where the Laplacian spectrum is often more effective than the adjacency spectrum in distinguishing graphs. 
In fact, the transition $q$-spectrum \jb{can} distinguish all the graphs with up to seven vertices, \jb{whereas the $q$-spectrum cannot}. 
% This gap between theory and computation poses a compelling open problem for future research.

{ On the other hand, there are non-isomorphic graphs that cannot be distinguished even by the transition $q$-spectrum. We show that, for strongly regular graphs, the magnitude, $q$-spectrum, and transition $q$-spectrum are all determined by the four parameters $(n,k,\lambda,\mu)$. Consequently, non-isomorphic strongly regular graphs with the same parameters provide examples with identical values of all three invariants\footnote{The first such example, a pair of the rook graph and the Shrikhande graph was found by Fable 5.}. Since infinitely many such pairs exist, this also yields infinitely many pairs of non-isomorphic graphs with identical transition $q$-spectra.
}

\smallskip
%!!!!!!!!!!!!!!!!!!!!!!!!!!!!!!!!!!!!!%%%%%%%%%%%%%%%%%%%%%%%%%%%%%%%%%%%%
\section{Preliminaries}\label{section_main_thm}
%!!!!!!!!!!!!!!!!!!!!!!!!!!!!!!!!!!!!!%%%%%%%%%%%%%%%%%%%%%%%%%%%%%%%%%%%%
%
First we give definitions of our invariants. 
Let $(X,d)$ be a finite metric space with $X=\{x_1,\dots,x_n\}$ and put $d_{ij}=d(x_i,x_j)$. 
The {\em similarity matrix} is given by $Z(q)=\left(q^{d_{ij}}\right)_{i,j}$. 
\begin{definition} \rm 
The {\em $q$-spectrum} of $X$ is a multiset (i.e. a set allowing multiple instances of elements) $[\la_1(q), \dots, \la_n(q)]$ of the eigenvalues of $Z(q)$. 
We express multisets by using brackets $[\dots]$ instead of $\{\dots\}$ hereafter. 
\end{definition}

We denote by $p(q;\la)$ the characteristic polynomial of $Z(q)$; $p(q;\la)=\left|\la I-Z(q)\right|$. 

\smallskip
Throughout this paper, graphs are assumed to be connected, undirected, and simple (without loops or multiple edges). By defining the distance between two vertices as the length of the shortest path connecting them, the vertex set of a graph forms a metric space. The adjacency matrix $A$ of a graph $G$ is a matrix whose $(i,j)$-entry is $1$ if vertices $v_i$ and $v_j$ are adjacent, and $0$ otherwise. The multiset of eigenvalues of this matrix are referred to as the {\em adjacency spectrum} of the graph, and graphs sharing the same adjacency spectrum are called {\em adjacency cospectral} graphs. 
Another fundamental matrix in spectral graph theory is the Laplacian matrix (see, for example, \cite{BH}). It is given by $L = D - A$, where $D$ is the degree matrix. The multiset of its eigenvalues is called the \emph{Laplacian spectrum}, and graphs sharing the same Laplacian spectrum are referred to as \emph{Laplacian cospectral} graphs.

Since for a graph $G$, 
\[
Z(q)=I+q\,A+h.o.t.,
\]
we have 
\begin{proposition}\label{q-spectrum_implies_spectrum}
The adjacency spectrum can be recovered from the $q$-spectrum by taking the limit of $(\lambda-1)/q$ as $q\to0$. 
\end{proposition}
Thus, the $q$-spectrum serves as a natural extension of the graph spectrum to general finite metric spaces. 

\smallskip
We also employ spectral invariants derived from the transition (stochastic) matrix associated with the similarity matrix. 

\begin{definition} \rm 
Let $X$ be a finite metric space and $Z(q)=\left(q^{d_{ij}}\right)_{i,j}$ be the similarity matrix as before. 
Let $\pi_i$ be the sum of the entries of the $i$-th row, $\pi_i=\sum_{j=1}^n q^{d_{ij}}$. 
Define the {\em row-stochastic transition  similarity matrix} by 
$P(q)=\left(q^{d_{ij}}/\pi_i\right)_{i,j}$. 
We call the eigenvalues of $P(q)$ the {\em transition (or stochastic) $q$-spectrum} of $X$. 
\end{definition}

We denote by $ps(q;\la)$ the characteristic polynomial of $P(q)$. 

As in spectral graph theory, the transition matrix is similar to a symmetric matrix, and hence diagonalizable with real eigenvalues in $[-1,1]$ (\cite{C,LPW}). 
To be precise, as $\pi_i\cdot{(P)}_{i,j}=q^{d_{ij}}=q^{d_{ji}}=\pi_j\cdot{(P)}_{j,i}$, the transition matrix $P$ is reversible. 
Put $\Pi=\mbox{\rm diag}\,(\pi_1,\dots,\pi_n)$, then $\Pi^{1/2}P\Pi^{-1/2}=\left(q^{d_{ij}}/\sqrt{\pi_i\pi_j}\right)_{i,j}$ is real symmetric, and therefore $P$ is diagonalizable and the eigenvalues (i.e. transition $q$-spectra) are all real numbers belonging to $[-1,1]$. 

As in the theory of random walks on graphs, 
the positivity of the transition similarity matrix implies that it is primitive (in particular, irreducible and aperiodic). 
Hence, by the Perron-Frobenius theorem (see, for example, \cite{Lo}), 
the largest eigenvalue of $P(q)$ is $1$, it has multiplicity one, 
and the corresponding eigenvector can be chosen to have strictly positive entries. 
Moreover, the powers of $P(q)$ converge to the stationary distribution as left eigenvector: 
\[
\lim_{k\to\infty} (P(q)^k)_{i,j}=v_j,
\]
where $(v_j)_j$ is the normalized Perron-Frobenius eigenvector.

\smallskip
Since for a graph $G$, 
\[
\lim_{q\to0}\frac{P-I}q=-D+A=-L,
\]
we have 
\begin{proposition}\label{transition_q-spectrum_L_spectrum_limit}
For graphs, the transition $q$-spectrum recovers the Laplacian spectrum by taking the limit of $-(\lambda-1)/q$ as $q\to0$. 
\end{proposition}
Thus, the transition $q$-spectrum serves as a natural extension of the Laplacian spectrum to general finite metric spaces.

\smallskip
Next we present the moduli space of unordered $n$-point metric spaces introduced by Roff and Yoshinaga (\cite{RY} Section 2). 
Put $N_2={n\choose2}$ and 
\[\mathcal{L}_n=\left\{(\de_{12},\de_{13},\dots,\de_{n-1\,n})\in (\RR_{>0})^{N_2}\,\left|\,
\begin{array}{l}
\,\de_{ij}+\de_{jk}\ge \de_{ik} \>\, \forall i,j,k\,(1\le i,j,k\le n),  \\[0.5mm]
\mbox{where we put }\>\de_{ii}=0,\,  \de_{ij}=\de_{ji}\>\,\forall i,j
\end{array}
\right.\right\}.\]
The symmetric group $\frak{S}_n$ acts on $\mathcal{L}_n$ by $\sigma\cdot(\de_{ij})=(\de_{\si(i)\si(j)})$. 
Then the {\em moduli space} is given by $\mathcal{M}_n=\mathcal{L}_n/\frak{S}_n$ equipped with the quotient topology. 

Finally, we introduce two kinds of genericity for metric spaces. 
\begin{definition} \rm 
We say that $X$ is {\em $3$-generic} if the sums of any three distinct pairwise distances are mutually distinct. 
\end{definition}
We remark the above condition of $3$-genericity is weaker than that of $p$-genericity with $p=3$ in Definition 2.3 of \cite{O24}. 

If $X$ is $3$-generic then the all pairwise distances are distinct. Therefore, if we put 
\[
\mathcal{L}_n^{3\mbox{-{\footnotesize gen}}}=\left\{(d_{12},d_{13},\dots,d_{n-1\,n})\in\mathcal{L}_n\,|\,X\>\mbox{ is $3$-generic}\right\}
\]
then the symmetric group $\frak{S}_n$ acts on $\mathcal{L}_n^{3\mbox{-{\footnotesize gen}}}$ freely. 
Moreover, the moduli space of $3$-generic $n$-point metric space, $\mathcal{M}_n^{3\mbox{-{\footnotesize gen}}}=\mathcal{L}_n^{3\mbox{-{\footnotesize gen}}}/\frak{S}_n$, is open dense in $\mathcal{M}_n$ since $\mathcal{L}_n\setminus\mathcal{L}_n^{3\mbox{-{\footnotesize gen}}}$ consists of the solutions of only finitely many linear equations. 

\begin{definition}\label{def_Q-generic} \rm(\cite{O24}) 
We say that $X$ is {\em rationally generic}, written $\QQ$-generic, if $d_{ij}$ $(i<j)$ are rationally independent. 
\end{definition}

Rationally generic spaces are generic in the sense that the set of rationally generic $n$-point metric space $\mathcal{M}_n^{\mbox{{\footnotesize $\QQ$-gen}}}$ has a measure-zero complement in $\mathcal{M}_n$, which is because the complement is a union of countably many codimension one subspaces. However it is not open dense in $\mathcal{M}_n$ since the complement includes $(\mathcal{L}_n\cap\QQ^{N_2})/\frak{S}_n$. 
Obviously rationally generic space is $3$-generic.

%!!!!!!!!!!!!!!!!!!!!!!!!!!!!!!!!!!!!!%%%%%%%%%%%%%%%%%%%%%%%%%%%%%%%%%%%%
\section{The $q$-spectrum}\label{section_main_thm}
%!!!!!!!!!!!!!!!!!!!!!!!!!!!!!!!!!!!!!%%%%%%%%%%%%%%%%%%%%%%%%%%%%%%%%%%%%

%!!!!!!!!!!!!!!!!!!!!!!!!!!!!!!!!!!!!!%%%%%%%%%%%%%%%%%%%%%%%%%%%%%%%%%%%%
\subsection{Identification of finite metric spaces}%\label{}
%!!!!!!!!!!!!!!!!!!!!!!!!!!!!!!!!!!!!!%%%%%%%%%%%%%%%%%%%%%%%%%%%%%%%%%%%%
%
We introduce results on distinguishing spaces by $q$-spectrum. 

First, for graphs, Proposition \ref{q-spectrum_implies_spectrum} implies 
\begin{proposition}\label{q-spectrum>spectrum}
The $q$-spectrum is strictly stronger than the adjacency spectrum for distinguishing graphs. 
\end{proposition}
The strictness follows from the existence of adjacency cospectral graphs with different $q$-spectra shown in Section \ref{section_experiments} (Table \ref{table}). 

Next, the $q$-spectrum becomes a complete invariant for generic finite metric spaces, to be precise, 
\begin{theorem}\label{main_thm}
Any $3$-generic metric space can be reconstructed from its $q$-spectrum. 
\end{theorem}

\smallskip
Finally, restricted to metric spaces on at most four points, $q$-spectrum is a complete invariant, namely, 
\begin{theorem}\label{thm_four-point} 
Any metric space on at most four points can be reconstructed from its $q$-spectrum. 
\end{theorem}
%

%!!!!!!!!!!!!!!!!!!!!!!!!!!!!!!!!!!!!!%%%%%%%%%%%%%%%%%%%%%%%%%%%%%%%%%%%%
\subsection{Combinatorial interpretation of the $q$-spectrum}\label{section_proof}
%!!!!!!!!!!!!!!!!!!!!!!!!!!!!!!!!!!!!!%%%%%%%%%%%%%%%%%%%%%%%%%%%%%%%%%%%%
%
The key observation is that the coefficients of the characteristic polynomial admit a combinatorial interpretation in terms of cycles of the metric space.

A sequence $\Ga=(x_{i_0},\dots,x_{i_k})$ of $k+1$ points of $X$ such that adjacent points are different is called a {\em $k$-path} of $X$. 
Its length is defined by the sum of $k$ pairwise distances; $\ell(\Ga)=\sum_{j=1}^k d_{i_{j-1}i_{j}}$. 

\begin{definition}\label{def_cycle} \rm 
Let $k\ge2$. 
By a {\em single $k$-cycle} we mean the equivalence class of sequences of $k$ mutually distinct points $\Ga=[x_{i_1},\dots,x_{i_k}]=(x_{i_1},\dots,x_{i_k})/\sim$, where the equivalence is generated by the cyclic shift;
\[
(x_{i_1},x_{i_2},\dots,x_{i_k})\sim(x_{i_2},\dots,x_{i_k},x_{i_1}). 
\]
The length of a single $k$-cycle is defined to be the sum of $k$ pairwise distances; 
$\ell(\Ga)=\sum_{j=1}^{k-1} d_{i_{j}i_{j+1}}+d_{i_ki_1}.$

By a {\em $k$-cycle} we mean a disjoint union of single cycles such that the sum of the numbers of vertices is equal to $k$. Its length is defined to be the sum of the lengths of constituent single cycles. 
Let $c(\Ga)$ be the number of connected components of a $k$-cycle $\Ga$. 
Let $\mathcal{C}_k(X)$ be the set of $k$-cycles of $X$. 
Put $\mathcal{C}_1(X)=\emptyset$. 
\end{definition}

Remark that both $k$-paths and $k$-cycles are oriented.

Let $\bar p(q;\mu)$ be a polynomial obtained from $p(q;\la)$ by setting $\mu=\la-1$; $\bar p(q;\mu)=p(q;\mu+1)$. 
Put $\mu_i(q)=\la_i(q)-1$, and let $\tau_k(q)$ $(1\le k \le n)$ be the $k$-th elementary symmetric polynomial of $\mu_1(q),\dots, \mu_n(q)$. 
{
}
Now we state the key lemma: 

\begin{lemma}\label{lem_tau}
There holds
\begin{equation}\label{tau_L}
\tau_k(q)=\sum_{\Ga\in\mathcal{C}_k(X)}{(-1)}^{c(\Ga)+k}q^{\,\ell(\Ga)} \qquad (1\le k\le n). 
\end{equation}
In particular, 
\begin{equation}\label{tau2_tau3}
\tau_2(q)=-\sum_{i<j}q^{2d_{ij}}, \quad \tau_3(q)=2\sum_{i<j<k} q^{d_{ij}+d_{jk}+d_{ki}}.
\end{equation}
\end{lemma}

\begin{proof}
When $k=1$, $\tau_1(q)={\rm trace}\,(I-Z(q))=0$. 

Suppose $k\ge2$. As there holds 
\begin{equation}\label{bar_p}
\bar p(q;\mu)=
\prod_{i=1}^n(\mu-\mu_i(q))
=\mu^n+\sum_{l=1}^n(-1)^l\tau_l(q)\,\mu^{n-l}
=\left|\begin{array}{cccc}
\mu & -q^{d_{12}} & \cdots & -q^{d_{1n}} \\
-q^{d_{21}} & \mu & & \vdots\\
\vdots & & \ddots & -q^{d_{n-1\,n}} \\[1mm]
-q^{d_{n1}} & \cdots &  -q^{d_{n\,n-1}} & \mu 
\end{array}
\right|, 
\end{equation}
the coefficient of $\mu^{n-k}$ is given by 
\begin{equation}\label{coeff_mu_n-k}
\sum_{1\le i_1<\dots<i_k\le n\>{}}\sum_{{}\>\{j_1,\dots,j_k\}=\{i_1,\dots,i_k\}, j_m\ne i_m (1\le m\le k)}
{\rm sgn}\left(\begin{array}{ccc}
i_1 & \cdots  & i_k \\
j_1 & \cdots  & j_k
\end{array}
\right)\left(-q^{d_{i_1j_1}}\right)\dots\left(-q^{d_{i_kj_k}}\right). 
\end{equation}

A permutation of $k$ letters such that no letters remain unchanged can be expressed as the product of several cyclic permutations of disjoint orbits with $k_i$ letters, where $\sum_i k_i=k$, and the expression is unique up to the order of the cycles. 
A cyclic permutation of $m$ letters $(i_1, \dots, i_m)$ corresponds to a {\sl single} $m$-cycle $\Ga_{i_1\dots i_m}=[x_{i_1},\dots,x_{i_m}]$. 
Since the signature of a cyclic permutation of $m$ letters is $(-1)^{m-1}$, we have 
\[
\displaystyle 
{\rm sgn}(i_1, \dots, i_m)
\big(-q^{d_{i_1i_2}}\big)\dots\big(-q^{d_{i_{m-1}i_m}}\big)
\big(-q^{d_{i_mi_1}}\big)
=-q^{\ell(\Ga_{i_1\dots i_m})}. 
\]
Since a permutation of $k$ letters corresponds to a $k$-cycle in the sense of Definition \ref{def_cycle}, the right hand side of \eqref{coeff_mu_n-k} is equal to 
\[
\sum_{\Ga\in\mathcal{C}_k(X)}{(-1)}^{c(\Ga)}q^{\,\ell(\Ga)},
\]
which, together with \eqref{coeff_mu_n-k}, implies \eqref{tau_L}. 

\smallskip
We remark that the coefficient $1/2$ of the right hand side of the second equality of \eqref{tau2_tau3} comes from the fact that cycles are oriented and therefore every triangle is counted twice. 
\end{proof}

{%
\begin{remark}\rm 
For $i\in\mathbb N$ and $\ell>0$, define
\[
C_{i;k}^{\ell}
=
\mathbb Z
\left\{
\Gamma\in\mathcal C_k(X)
\;:\;
c(\Gamma)=i,\;
\ell(\Gamma)=\ell
\right\}.
\]
We regard $i$ as the homological degree and $\ell$ as the internal degree.
The graded Euler characteristic of $\displaystyle C_k=\{C_{i;k}^{\ell}\}_{i,\ell}$
is defined by
\[
\chi_q(C_k)
=
\sum_{i\in\mathbb N,\ \ell>0} (-1)^i \operatorname{rank} C_{i;k}^{\ell}\, q^{\ell}.
\]
By the lemma above, 
\[
\tau_k(q)=(-1)^k\chi_q(C_k).
\]
Thus, the elementary symmetric polynomials of the shifted $q$-spectrum $\la_i(q)-1$ admit a categorified interpretation as graded Euler characteristics of cycle spaces. 
\end{remark}
}

\begin{corollary}\label{cor_2_3-pt}
Complete graphs, two-point metric spaces and three-point metric spaces can be identified by the $q$-spectrum. 
\end{corollary}

\begin{corollary}\label{cor_diameter}
Spaces with different diameters have different $q$-spectra. 
\end{corollary}

This is because the largest power appearing in $\tau_2(q)$ is the diameter of $X$. 

\begin{proposition}\label{lem_multisets}
The $q$-spectrum determines not only the multiset of pairwise distances $\mathcal{S}_1=[d_{ij}]_{i<j}$ but also the multiset of triangle perimeters $\mathcal{S}_3^\triangle=[d_{ij}+d_{jk}+d_{ki}]_{i<j<k}$. 
\end{proposition}

\begin{proof}
Define $\tilde \tau_k(t)$ by $\tilde \tau_k(t)=\tau_k(e^{-t})$  $(k=2,3)$. Then 
\[
\sum_{i<j}e^{-2d_{ij}t}=-\tilde \tau_2(t), \quad 
\sum_{i<j<k} e^{-(d_{ij}+d_{jk}+d_{ki})t}=\frac12 \tilde \tau_3(t), 
\]
which implies 
\[
\sum_{i<j} \, {(d_{ij})}^m=\frac{-1}{(-2)^m}\,\tilde\tau_2^{(m)}(0), \quad
\sum_{i<j<k}\,{(d_{ij}+d_{jk}+d_{ki})}^m=\frac{(-1)^{m}}2\,\tilde\tau_3^{(m)}(0). 
\]
The left hand sides above are the Newton polynomials of $d_{ij}$ and $d_{ij}+d_{jk}+d_{ki}$ respectively. 
Since the elementary symmetric polynomials can be obtained from the Newton polynomials, we obtain $\mathcal{S}_1=[d_{ij}]_{i<j}$ from $\tilde \tau_2^{(m)}(0)$ $(1\le m\le {n\choose2})$ and $\mathcal{S}_3^\triangle=[d_{ij}+d_{jk}+d_{ki}]_{i<j<k}$ from $\tilde \tau_3^{(m)}(0)$ $(1\le m\le {n\choose3})$. 
\end{proof}

The reason the condition required for the $q$-spectrum in the space distinguishing theorem is weaker than that for %the magnitude or 
the transition $q$-spectrum (Theorem \ref{thm_Q-gen_transition}) is that $\mathcal{S}_3^\triangle$ can be recovered without assuming genericity.

%!!!!!!!!!!!!!!!!!!!!!!!!!!!!!!!!!!!!!%%%%%%%%%%%%%%%%%%%%%%%%%%%%%%%%%%%%
\subsection{Proof for $3$-Generic case}\label{subsection_generic}
%!!!!!!!!!!!!!!!!!!!!!!!!!!!!!!!!!!!!!%%%%%%%%%%%%%%%%%%%%%%%%%%%%%%%%%%%%
%

\bigskip\noindent
{\bfseries Proof of Theorem {\rm \bf \ref{main_thm}}.} 
Propostion \ref{lem_multisets} implies that $\mathcal{S}_1$ and $\mathcal{S}_3^\triangle$ can be obtained from the $q$-spectrum. 
We show inductively that a $3$-generic finite metric space can be constructed from $\mathcal{S}_1$ and $\mathcal{S}_3^\triangle$. 
Put $\mathcal{S}_1=[a_i]_{i=1}^{N_2}$, where $N_2={n\choose2}$, and $\mathcal{S}_3^\triangle=[t_j]_{j=1}^{N_3}$, where $N_3={n\choose3}$. 
Note that the $3$-genericity condition implies that both multisets $\mathcal{S}_1$ and $\mathcal{S}_3^\triangle$ become sets, and that for each element $t_m\in \mathcal{S}_3^\triangle$ there is a unique triple $\{a_i,a_j,a_k\}\subset \mathcal{S}_1$ that satisfies $t_m=a_i+a_j+a_k$. 

(i) Take any element of $\mathcal{S}_3^\triangle$, say $t_1$. 
Let $\{a,b,c\}\subset S_1$ be the unique triple that satisfies $t_1=a+b+c$. 
Name the vertices $A,B$ and $C$ in the standard way, and choose the edge $BC$ with length $a$ as the base. 
Put $\mathcal{S}_3^\triangle(a)=\{t\in \mathcal{S}_3^\triangle\,|\,t=a+x+y, \, x,y\in \mathcal{S}_1\}$. 
We may assume, by renumbering the subscripts of $t_j$ if necessary, that $\mathcal{S}_3^\triangle(a)=\{t_1, \dots, t_{n-2}\}$. 

\begin{figure}[htbp]
\begin{center}
\includegraphics[width=.45\linewidth]{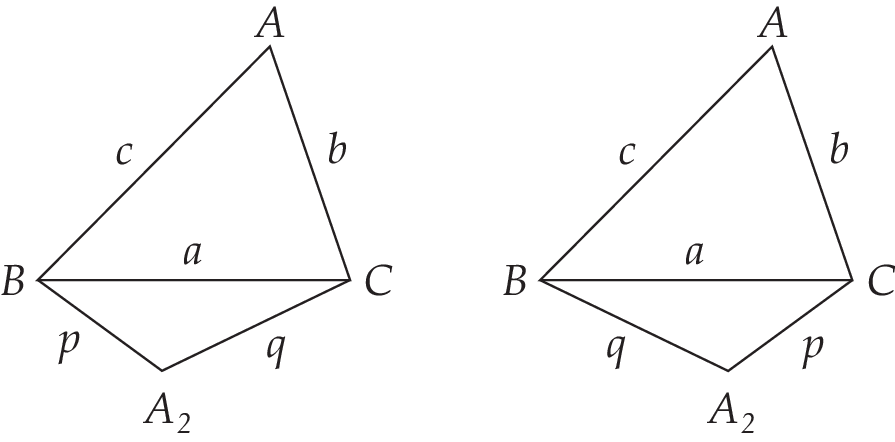}
\caption{}
\label{attach_a_triangle}
\end{center}
\end{figure}

(ii) Take $t_2$. Suppose $t_2=a+p+q$, where $p,q\in \mathcal{S}_1\setminus\{a,b,c\}$. 
There are two ways to attach a triangle with pairwise distances $a,p,q$, which we denote by $\triangle (a,p,q)$, to the edge $BC$. 
We choose one of them by attaching the edge with length $p$ to the vertex $B$ if there is $t\in \mathcal{S}_3^\triangle\setminus\{t_1,t_2\}$ such that $t=p+c+x$ for some $x\in \mathcal{S}_1\setminus\{a,b,c,p,q\}$, and to $C$ otherwise (Figure \ref{attach_a_triangle}). 
Put $A_2$ to be the vertex of  $\triangle (a,p,q)$ which is different from $B$ and $C$, then $AA_2=x$. 

(iii) Take $t_3\in \mathcal{S}_3^\triangle(a)$. Assume $t_3=a+r+s$, where $r,s\in \mathcal{S}_1\setminus\{a,b,c,p,q\}$. 
We can determine the way how $\triangle (a,r,s)$ is attached to the edge $BC$ in the same way as above. 
Let $A_3$ be the vertex of $\triangle (a,r,s)$ which is different from $B$ and $C$. 
Then the length $AA_3$ is given by the unique element $u\in S_1\setminus\{a,b,c,p,q,r,s\}$ that satisfies $c+r+u\in \mathcal{S}_3^\triangle$ or $b+r+u\in \mathcal{S}_3^\triangle$ depending on the way how $\triangle (a,r,s)$ is attached to $BC$. 
The length $A_2A_3$ can be obtained in the same way. 

(iv) Repeat the above procedure. 
{\hfill{\small $\square$}\par\medskip}
\medskip

We remark that the theorem does not exclude the possibility that a $3$-generic space and a non $3$-generic space have the same $q$-spectrum. 
We also remark that even without the assumption of genericity, the multiset of pairwise distances $\mathcal{S}_1$ is determined by the $q$-spectrum and hence that the number of isometry classes of spaces with a given $q$-spectrum is bounded universally.

%!!!!!!!!!!!!!!!!!!!!!!!!!!!!!!!!!!!!!%%%%%%%%%%%%%%%%%%%%%%%%%%%%%%%%%%%%
\subsection{Proof for Four-point sets case}\label{subsection_four_pts}
%!!!!!!!!!!!!!!!!!!!!!!!!!!!!!!!!!!!!!%%%%%%%%%%%%%%%%%%%%%%%%%%%%%%%%%%%%
%

\noindent
{\bfseries Proof of Theorem {\rm \bf \ref{thm_four-point}}.} 
The statement for the case when the number of points is two or three follows from Proposition \ref{lem_multisets} as was stated in Corollary \ref{cor_2_3-pt}. 
We give proof for four-point sets in what follows. 

We use the $k$-th elementary symmetric polynomials $\tau_k(q)$ of $\mu_k(q)$ $(k=2,3,4)$, where $\mu_i(q)=\la_i(q)-1$. 
By Proposition \ref{lem_multisets}, from $\tau_2(q)$ and $\tau_3(q)$ we obtain a multiset of pairwise distances and a multiset of perimeters of four triangles, which we denote by $\mathcal{S}_1=[a,b,c,d,f,g]$ and $\mathcal{S}_3^\triangle=[t_1,t_2,t_3,t_4]$ respectively. 
Put $s$ to be the sum of six pairwise distances $s=a+b+c+d+f+g$. The proof consists of two steps. 

\begin{figure}[htbp]
\begin{center}
\includegraphics[width=.55\linewidth]{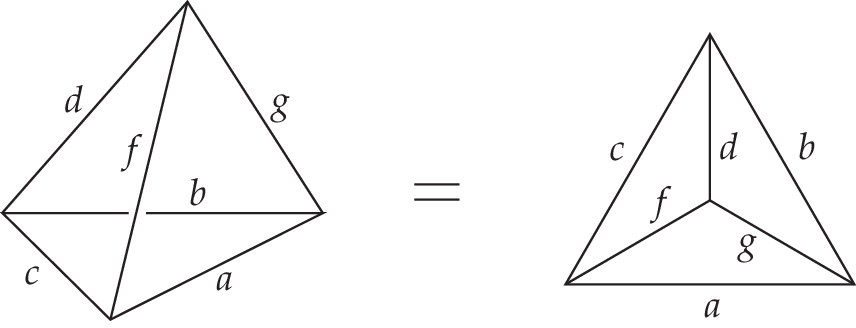}
\caption{}
\label{tetrahedron}
\end{center}
\end{figure}

\smallskip
Step 1. 
First we show that we can obtain a multiset of the sums of the lengths of {\sl opposite} (i.e. non-adjacent) edges, which we denote by $\mathcal{S}_{\rm opp}=[\a,\beta,\ga]$. 
Note that $\a+\beta+\ga=s$. 
We label the pairwise distances $a,\dots,g$ as shown in Figure \ref{tetrahedron}, and 
put 
\begin{equation}\label{alphabetagamma}
\a=a+d, \>\> \beta=b+f, \>\> \ga=c+g.
\end{equation}
Then the lengths of single $4$-cycles are given by 
\[\a+\beta=a+b+d+f, \>\> \beta+\ga=b+c+f+g, \>\> \a+\ga=a+c+d+f,\]
and the lengths of $4$-cycles that consist of two single $2$-cycles are given by 
$2\a, \> 2\beta$ and $2\ga$. 
Therefore $\tau_4(q)$, the constant term of $\bar p(q;\mu)$, is given by 
\begin{equation}\label{tau_4}
-2q^{\a+\beta}-2q^{\beta+\ga}-2q^{\a+\ga}+q^{2\a }+q^{2\beta }+q^{2\ga }. 
\end{equation}
There are at most four cases (i) - (iv) depending on whether none, exactly two, or three of $\alpha, \beta$ and $\gamma$ conincide. 

(1-i) $\tau_4(q)$ is of the form $-2\sum_{i=1}^3q^{A_i}+\sum_{i=j}^3q^{B_j},$ 
where $A_i, B_j$ are mutually distinct. 
This case can occur if and only if $\a,\beta$ and $\ga$ are mutually distinct and they do not form an arithmetic progression. 
In this case $\mathcal{S}_{\rm opp}$ is given by $\mathcal{S}_{\rm opp}=[B_1/2, B_2/2, B_3/2]$. 

(1-ii) $\tau_4(q)$ is of the form 
\[-2q^{A_1'}-2q^{A_2'}-q^{A_3'}+q^{B_1'}+q^{B_2'},\]
where $A_1',A_2',A_3',B_1'$ and $B_2'$ are mutually distinct. 
This case can occur if and only if $\a,\beta$ and $\ga$ are mutually distinct and they form an arithmetic progression. 
We may assume without loss of generality that $\a+\beta=2\ga$. Then cancellation between $-q^{\a+\beta}$ and $q^{2\ga }$ occurs in \eqref{tau_4}, and $\mathcal{S}_{\rm opp}$ is given by $\mathcal{S}_{\rm opp}=[B_1'/2, B_2'/2, s-(B_1'+B_2')/2]$. 

(1-iii) $\tau_4(q)$ is of the form $-4q^{A}+q^{B},$ where $A$ and $B$ are distinct. 
This case can occur if and only if exactly two of $\a,\beta$ and $\ga$ coincide. 
We may assume without loss of generality that $\a=\beta\ne\ga$. 
Then cancellation between $-2q^{\a+\beta}$ and $q^{2\a }+q^{2\beta }$ occurs in \eqref{tau_4}, 
and $\mathcal{S}_{\rm opp}$ is given by $\mathcal{S}_{\rm opp}=[B/2, (s-B/2)/2, (s-B/2)/2]$. 

(1-iv) $\tau_4(q)$ is of the form $-3q^{A'}$. 
This case can occur if and only if $\a,\beta$ and $\ga$ are all equal. 
In this case $\mathcal{S}_{\rm opp}$ is given by $\mathcal{S}_{\rm opp}=[A'/2, A'/2, A'/2]$. 

\smallskip
Step 2. We show that if $X$ and $X'$ both realize the same $\mathcal{S}_1, \mathcal{S}_3^\triangle$ and $\mathcal{S}_{\rm opp}$ then $X$ and $X'$ are isometric. 
This is done by listing up all the possible configurations (i.e. the ways to assign elements of $\mathcal{S}_1$ to edges of a tetrahedron) with given data of $\mathcal{S}_1$ and $\mathcal{S}_{\rm opp}$ and then by showing that any two of them are isometric if they have the same $\mathcal{S}_3^\triangle$. 
The remaining part of the proof consists of an exhaustive case analysis. 
The details are deferred to Appendix. 
{\hfill{\small $\square$}\par\medskip}
\medskip

%!!!!!!!!!!!!!!!!!!!!!!!!!!!!!!!!!!!!!%%%%%%%%%%%%%%%%%%%%%%%%%%%%%%%%%%%%
\section{Transition $q$-spectrum}
%!!!!!!!!!!!!!!!!!!!!!!!!!!!!!!!!!!!!!%%%%%%%%%%%%%%%%%%%%%%%%%%%%%%%%%%%%

Recall that the transition $q$-spectrum is given by the eigenvalues of $P(q)=\left(q^{d_{ij}}/\pi_i\right)_{i,j}$, where $\pi$ is the sum of the entries of the $i$-th row, $\pi_i=\sum_{j=1}^n q^{d_{ij}}$. 

Analogous results to the theorems on distinguishing spaces by the $q$-spectrum hold in a weaker form.
Proposition \ref{transition_q-spectrum_L_spectrum_limit} implies 
\begin{proposition}\label{transition_q-spectrum>L_spectrum}
The transition $q$-spectrum is strictly stronger than the Laplacian spectrum for distinguishing graphs. 
\end{proposition}
The strictness follows from the existence of graphs with the same Laplacian spectrum and different transition $q$-spectra shown in Section \ref{section_experiments} (Table \ref{table}). 

Transition $q$ spectrum encodes combinatorial data of metric spaces as follows. 
The trace of $P(q)$ is equal to the sum of the transition $q$-spectra. 
Put $f(q)={\rm tr}\,P(q)$, then it can be expanded in a generalized polynomial of $q$ using 
\[
s_i=\pi_i-1=\sum_{j\ne i}q^{d_{ij}},
\]
as 
\begin{eqnarray}
f(q)={\rm tr}\,P(q) &=& \displaystyle \sum_i\frac1{1+s_i}
=n-\sum_is_i+\sum_i s_i^2-\sum_i s_i^3 +\dots. \label{ser_exp_trace}
\end{eqnarray}
\begin{theorem}\label{thm_Q-gen_transition} 
Any rationally generic metric space (Definition \ref{def_Q-generic}) can be reconstructed from its transition $q$-spectrum. 
\end{theorem}

\begin{proof}
Let $\mathcal{P}$ denote the multiset of powers of $q$ appearing in equation \eqref{ser_exp_trace}. 
Rational independence enshures that no two terms appearing in \eqref{ser_exp_trace} can coincide. 
Since 
\[
-\sum_i s_i = -2\sum_{i<j} q^{d_{ij}},
\]
and no other term in \eqref{ser_exp_trace} can have $-2$ as its coefficient, the set of pairwise distances $\mathcal{S}_1$ can be obtained as the set of exponents of the terms having coefficient $-2$ in \eqref{ser_exp_trace}. 

Next, we determine the incidence relation of the ``edges''. Among the sums of two distinct elements of $\mathcal{S}_1$, those appearing in $\mathcal{P}$ give rise, through the term $\sum_i s_i^2$, to the set  
\[
\hat{\mathcal{S}}_2=\{\,d_{ij}+d_{ik}\mid i\neq j\neq k\neq i\,\}. 
\]
Suppose that
\[
a+b,\quad b+c,\quad c+a
\]
belong to $\hat{\mathcal{S}}_2$,
where $a,b,c$ are distinct elements of $\mathcal{S}_1$. 
Then the ``edges'' of lengths $a,b,c$ form a triangle. 

Now consider another triangle $a,b',c'$ sharing the edge $a$, where $b',c'\notin\{b,c\}.$ 
There are two possible ways to glue the triangle 
$a,b',c'$ to the triangle $a,b,c$ along the edge $a$. 
We declare that the edges $c$ and $c'$ are incident to the same vertex 
if and only if $c+c'\in \hat{\mathcal{S}}_2.$ 
This determines the gluing uniquely. 

Repeating this procedure determines all incidences between vertices and edges, 
and therefore reconstructs the metric space. 
\end{proof}

\begin{theorem}\label{thm_transition_3_pts}
Any metric space on at most three points can be reconstructed from its transition $q$-spectrum. 
\end{theorem}

\begin{proof}
Put $[a,b,c]=[d_{12},d_{13},d_{23}]$ for simplicity's sake. From 
\[
\begin{array}{rcl}
f'(1)&=&\displaystyle -\frac29(a+b+c), \\[4mm]
f''(1)&=&\displaystyle \frac2{27}\left(3(a+b+c)-(a^2+b^2+c^2)+2(ab+bc+ca)\right), \\[4mm]
f'''(1)&=&\displaystyle \frac2{27}\left(-6(a+b+c)+3(a^2+b^2+c^2)-6(ab+bc+ca)+(a^3+b^3+c^3)\right),
\end{array}
\]
one can obtain the elementary symmetric polynomials of $a,b$ and $c$. 
\end{proof}

We remark that this condition of Theorems \ref{thm_Q-gen_transition} is the same as the one used in the author's previous result on magnitude \cite{O24} (see Section \ref{section_experiments} for the definition of the magnitude). 
Note that the conditions of Theorems \ref{thm_Q-gen_transition} and \ref{thm_transition_3_pts} are more restrictive than that of Theorems \ref{main_thm} and \ref{thm_four-point}.

{{
%!!!!!!!!!!!!!!!!!!!!!!!!!!!!!!!!!!!!!%%%%%%%%%%%%%%%%%%%%%%%%%%%%%%%%%%%%
\section{Strongly regular graphs with identical similarity spectra}\label{section_srg}
%!!!!!!!!!!!!!!!!!!!!!!!!!!!!!!!!!!!!!%%%%%%%%%%%%%%%%%%%%%%%%%%%%%%%%%%%%

We now examine a limitation of the distinguishing power of the $q$-spectrum and the transition $q$-spectrum. We show that strongly regular graphs with the same parameters have identical $q$-spectra and transition $q$-spectra, and hence obtain infinitely many pairs of non-isomorphic graphs that cannot be distinguished by either invariant.

A {\em strongly regular graph}, denoted by $SRG(n,k,\lambda, \mu)$, is a graph such that 
\begin{itemize}
\item the number of vertices is equal to $n$,
\item the degree (the number of neighbours) is equal to $k$ for each vertex,
\item every two adjacent vertices have $\lambda$ common neighbours, 
\item every two non-adjacent vertices have $\mu$ common neighbours.
\end{itemize}

\begin{proposition}\label{prop_srg}
Let $G$ be a strongly regular graph with the parameters 
$(n,k,\lambda,\mu)$ and diameter $2$. 
Then its magnitude, $q$-spectrum and transition $q$-spectrum of $G$ are determined solely by $(n,k,\lambda,\mu)$. 
\end{proposition}

\begin{proof}
Since $G$ is regular with diameter $2$, every row sum of $Z(q)$ is equal to 
\begin{equation}\label{def_pi}
\pi(q)=1+kq+(n-k-1)q^2,
\end{equation}
and hence the formal magnitude is given by 
\begin{equation}\label{magnitude_4}
\textrm{Mag}(G)=\frac{n}{\pi(q)}=\frac{n}{1+kq+(n-k-1)q^2}  \nonumber
\end{equation}
since $(\pi(q))^{-1}\vect 1$, where $\vect 1={}^t(1,\dots,1)$, is a weighting (cf. Speyer, \cite{L13} Proposition 2.1.5). 

The transition similarity matrix is given by 
\begin{equation}\label{P-Z_4}
P(q)=\frac{1}{\pi(q)}Z(q),
\end{equation}
which implies that the transition $q$-spectrum are obtained from the $q$-spectrum by multiplying $(\pi(q))^{-1}$. 

It remains to express the $q$-spectrum in terms of $(n,k,\lambda,\mu)$. 
For this purpose, we first investigate the eigenvalues of the adjacency matrix $A$ and their multiplicities. 
Let $J$ be the matrix all the elements of which are equal to $1$. 
Since $G$ has diameter $2$, 
\begin{equation}\label{Z-A}
Z(q)=I+qA+q^2(J-I-A)
=(1-q^2)I+q(1-q)A+q^2J.
\end{equation}

Since $G$ is $k$-regular, it has an eigenvalue $k$ with an eigenvector $\vect 1$, and its multiplicity is $1$. 
The reason of the multiplicity $1$ is as follows. 
Since the Laplacian is given by $L=kI-A$, $Av=kv$ holds if and only if $Lv=0$. On the other hand, the dimension of the kernel of $L$ is equal to the number of connected components of the graph $G$ (Proposition 13.7 of \cite{BH}), which implies $\dim \ker (A-kI)=1$. 

Next, since $G$ has diameter $2$ 
\[
A^2=(k-\mu)I+(\lambda-\mu)A+\mu J,
\]
which implies that on $\vect 1^\perp$ there holds 
\[
A^2-(\lambda-\mu)A-(k-\mu)I=O.
\]
Therefore, if $t$ $(t\ne k)$ is an eigenvalue of $A$ then $t$ must satisfy 
\begin{equation}\label{quadratic_eq}
t^2-(\lambda-\mu)t-(k-\mu)=0.
\end{equation}
The discriminant 
\begin{equation}\label{discriminant}
\Delta=(\lambda-\mu)^2+4(k-\mu)
\end{equation}
is non-negative since $k\ge\mu$. 
Moreover, $\Delta=0$ implies $k=\mu=\lambda$, which cannot occur since $\lambda\le k-1$. 
Therefore $\Delta>0$. 

Let $t_1$ and $t_2$ be solutions to \eqref{quadratic_eq};
\begin{equation}\label{t1t2}
t_1,t_2=\frac{(\lambda-\mu)\pm\sqrt\Delta}2.
\end{equation}
They are eigenvalues of $A$. 
Let $m_1$ and $m_2$ be the multiplicity. Since
\[\begin{array}{rcl}
m_1+m_2+1&=&n,\\
t_1m_1+t_2m_2+k&=&\textrm{tr} A=0,
\end{array}\]
we have 
\begin{equation}\label{m1m2}
m_1,m_2=\frac12\left[(n-1)\mp\frac{2k+(n-1)(\lambda-\mu)}{\sqrt\Delta}\right].
\end{equation}
Now the adjacent spectrum of $G$ is given by $\textrm{Spec}_A(G)=\{k^1,t_1^{m_1},t_2^{m_2}\}.$ 

If $v$ is an eigenvector of $A$ that belongs to the eigenvalue $t_i$, $v\perp \vect 1$, i.e. $Jv=0$, therefore, \eqref{Z-A} implies 
\[
Z(q)v=\left((1-q^2)I+q(1-q)\,t_i\right)v.
\]
It follows that the $q$-spectrum and the transition $q$-spectrum of $G$ are expressed in terms of $n,k,\lambda,\mu$ as 
\[\begin{array}{rcl}
\textrm{Spec}_q(G)&=&\displaystyle \left\{\pi(q)^1, \> \left((1-q^2)+q(1-q)\,t_1\right)^{m_1}, \> \left((1-q^2)+q(1-q)\,t_2\right)^{m_2}\right\}, \\[3mm]
\textrm{Spec}_{{\rm tr.} q}(G)&=&\displaystyle \left\{1^1, \> \left(\frac{(1-q^2)+q(1-q)\,t_1}{\pi(q)}\right)^{m_1}, \> \left(\frac{(1-q^2)+q(1-q)\,t_2}{\pi(q)}\right)^{m_2}\right\},
\end{array}\]
where $\pi(q), t_1, t_2, m_1, m_2$ are given by \eqref{def_pi}, \eqref{discriminant}, \eqref{t1t2}, and \eqref{m1m2}. 
We remark that the second equality follows from the first by \eqref{P-Z_4}. 

\end{proof}

Since there are infinitely many pairs of non-isometric strongly regular graphs with the same parameter and with diameter $2$ (\cite{AH}), we have 

\begin{corollary}
There are infinitely many pairs of non-isomorphic graphs that share identical transition $q$-spectra, $q$-spectra, and magnitudes. 
\end{corollary}

We show the simplest example, a rook graph (Figure \ref{rook}) and a Shrikhande graph (Figure \ref{Shrikhande}), both of which have the parameters $(16,6,2,2)$. 

\begin{figure}[htbp]
\begin{center}
\begin{minipage}{.45\linewidth}
\begin{center}

\begin{tikzpicture}[
  scale=0.8,
  vtx/.style={circle,draw,fill=white,inner sep=0pt,minimum size=7mm,
              font=\footnotesize},
  rowe/.style={thick},
  cole/.style={semithick,dashed}
]
 
\foreach \i in {0,1,2,3}{
  \foreach \j in {0,1,2,3}{
    \node[vtx] (v\i\j) at ({2.4*\j},{-2.4*\i}) {\i\j};
  }
}
 
% ---- edges inside each row (solid) ----
\foreach \i in {0,1,2,3}{
  \draw[rowe] (v\i0)--(v\i1);
  \draw[rowe] (v\i1)--(v\i2);
  \draw[rowe] (v\i2)--(v\i3);
  \draw[rowe] (v\i0) to[bend left=28]  (v\i2);
  \draw[rowe] (v\i1) to[bend right=28] (v\i3);
  \draw[rowe] (v\i0) to[bend left=32]  (v\i3);
}
 
% ---- edges inside each column (dashed) ----
\foreach \j in {0,1,2,3}{
  \draw[cole] (v0\j)--(v1\j);
  \draw[cole] (v1\j)--(v2\j);
  \draw[cole] (v2\j)--(v3\j);
  \draw[cole] (v0\j) to[bend left=28]  (v2\j);
  \draw[cole] (v1\j) to[bend right=28] (v3\j);
  \draw[cole] (v0\j) to[bend left=32]  (v3\j);
}
 
\end{tikzpicture}

\caption{rook graph $K_4\square K_4$}
\label{rook}
\end{center}
\end{minipage}
\hskip 0.4cm
\begin{minipage}{.45\linewidth}
\begin{center}

\begin{tikzpicture}[
  scale=0.7,
  vtx/.style={circle,draw,fill=white,inner sep=0pt,minimum size=7mm,
              font=\footnotesize},
  gho/.style={circle,draw,densely dashed,fill=white,inner sep=0pt,
              minimum size=7mm,font=\footnotesize},
  hor/.style={thick},
  ver/.style={semithick,dashed},
  dia/.style={very thick,densely dotted}
]
 
% fundamental square of the torus
\draw[gray,very thin] (0,0) rectangle ({2.2*4},{2.2*4});
 
% ---- the 16 real vertices ----
\foreach \a in {0,1,2,3}{
  \foreach \b in {0,1,2,3}{
    \node[vtx] (u\a\b) at ({2.2*\a},{2.2*\b}) {\a\b};
  }
}
 
% ---- ghost copies: right border (a=4) and top border (b=4) ----
\foreach \b in {0,1,2,3}{
  \node[gho] (u4\b) at ({2.2*4},{2.2*\b}) {0\b};
}
\node[gho] (u44) at ({2.2*4},{2.2*4}) {00};
\foreach \a in {0,1,2,3}{
  \node[gho] (u\a4) at ({2.2*\a},{2.2*4}) {\a0};
}
 
% ---- generator (1,0): horizontal, solid ----
\foreach \a/\aa in {0/1,1/2,2/3,3/4}{
  \foreach \b in {0,1,2,3}{
    \draw[hor] (u\a\b)--(u\aa\b);
  }
}
 
% ---- generator (0,1): vertical, dashed ----
\foreach \b/\bb in {0/1,1/2,2/3,3/4}{
  \foreach \a in {0,1,2,3}{
    \draw[ver] (u\a\b)--(u\a\bb);
  }
}
 
% ---- generator (1,1): diagonal, dotted ----
\foreach \a/\aa in {0/1,1/2,2/3,3/4}{
  \foreach \b/\bb in {0/1,1/2,2/3,3/4}{
    \draw[dia] (u\a\b)--(u\aa\bb);
  }
}
 
\end{tikzpicture}

\caption{Shrikhande graph, drawn on a torus}
\label{Shrikhande}
\end{center}
\end{minipage}
\end{center}
\end{figure}

}}

%!!!!!!!!!!!!!!!!!!!!!!!!!!!!!!!!!!!!!%%%%%%%%%%%%%%%%%%%%%%%%%%%%%%%%%%%%
\section{Computational experiments}\label{section_experiments}
%!!!!!!!!!!!!!!!!!!!!!!!!!!!!!!!!!!!!!%%%%%%%%%%%%%%%%%%%%%%%%%%%%%%%%%%%%

We performed exhaustive computations, using Mathematica and Maple, in order to examine how effectively the $q$-spectrum and the transition $q$-spectrum distinguish non-generic finite metric spaces.

For comparison, we also consider magnitude, another invariant associated with the similarity matrix. 
The {\em formal magnitude} is defined to be the sum of all the entries of ${Z(q)}^{-1}$, where the inverse matrix is taken in the field of fractions of the integral domain consisting of ``{\em generalized polynomial}\,'' of $q$ that allow non-integer exponents (Leinster \cite{L19}). 

We also consider spectrum of the distance matrix $\left(d_{ij}\right)_{i,j}$, which we call the $d^{(1)}$-spectrum for graphs on vertices less than or equal to 7. 

We observed that \jb{among graphs with at most seven vertices, }
the transition $q$-spectrum distinguishes all the examples in our exhaustive computations, the transition $q$-spectrum is stronger than the $q$-spectrum, and the $q$-spectrum is stronger than the magnitude. 
Recall that, theoretically, the $q$-spectrum is stronger than the adjacent spectrum and the transition $q$-spectrum is stronger than the Laplacian spectrum (Propositions \ref{q-spectrum>spectrum} and \ref{transition_q-spectrum>L_spectrum}).
We also observed that neither the magnitude nor the adjacency spectrum is uniformly stronger than the other in terms of distinguishing power.

\medskip
Here we only show the results on graphs with vertices less than or equal to 7
The computation on $q$-spectrum, transition $q$-spectrum, $d^{(1)}$-spectrum, and magnitude was carried out using Mathematica by Tsuriba during his master's program (\cite{T}) (The reader is referred to \cite{AA,HS}, for example, for the results on the adjacent and Laplacian spectra). 
The built-in graph dataset in Mathematica for graphs with up to seven vertices was utilized.

\begin{table}[htbp]
 \begin{center}
\begin{tabular}{|c|r|r|r|r|r|r|r|} \hline
$n$ & $\#$ graphs & Adj. spec. & $q$-spec. & Lap. spec. & Trans. $q$ & Dist. $d^{(1)}$ & Mag. \\ \hline \hline
 4 & 6 & 0 & 0 & 0 & 0 & 0 & 2 \\ \hline
5 & 21 & 0 & 0 & 0 & 0 & 0 & 8 \\ \hline
6 & 112 & 2 & 0 & 4 & 0 & 0 & 63 \\ \hline
7 & 853 & 63 & 22 & 115 & 0 & 22 & 551 \\ \hline
 \end{tabular}
   \caption{Number of cospectral and co-magnitude graphs. 
The abbreviations in the top row of the table stand for the adjacency spectrum, $q$-spectrum, Laplacian spectrum, transition $q$-spectrum, $d^{(1)}$-spectrum, and magnitude respectively. }
\label{table}
 \end{center}
\end{table}

\begin{figure}[htbp]
\begin{center}
\includegraphics[width=.4\linewidth]{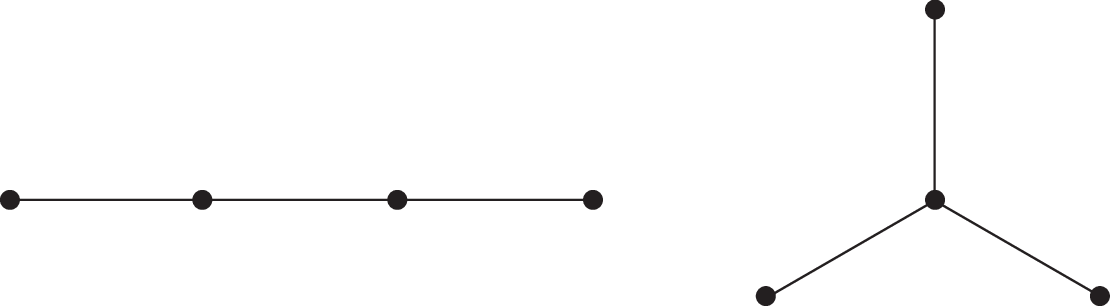}
\caption{Same magnitude with distinct spectra}
\label{two_graphs}
\end{center}
\end{figure}

There are six graphs on four vertices, and their graph spectra, both adjacent and Laplacian, are all distinct; consequently, both their $q$-spectra and transition $q$-spectra are also distinct. On the other hand, there is exactly one pair among them that shares the same magnitude (Figure \ref{two_graphs}).
It follows that neither the adjacency spectrum nor the Laplacian spectrum of graphs is weaker than magnitude in terms of distinguishing power. Consequently, the same is true for both the $q$-spectrum and the transition $q$-spectrum.

\begin{figure}[htbp]
\begin{center}
\begin{minipage}{.45\linewidth}
\begin{center}
\includegraphics[width=0.9\linewidth]{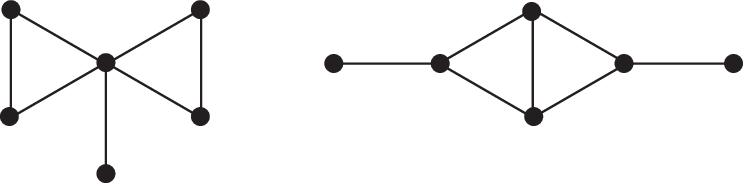}
\caption{Adjacensy cospectral, but not Laplacian cospectral} 
\label{cospectral-2}
\end{center}
\end{minipage}
\hskip 0.4cm
\begin{minipage}{.45\linewidth}
\begin{center}

\begin{tikzpicture}[scale=0.9,
  v/.style={circle, draw=black, fill=black, inner sep=1.3pt},
  e/.style={line width=0.8pt, draw=black}
]

% ----- G1 -----

\node[v] (a0) at (0,1) {};
\node[v] (a1) at (0,0) {};
\node[v] (a2) at (1,1) {};
\node[v] (a3) at (1,0) {};
\node[v] (a4) at (1,-1) {};
\node[v] (a5) at (2,0) {};

\draw (a0)--(a2);
\draw (a1)--(a2);
\draw (a1)--(a3);
\draw (a1)--(a4);
\draw (a2)--(a5);
\draw (a3)--(a5);
\draw (a4)--(a5);

% ----- G2 -----
\begin{scope}[xshift=4.2cm]

\node[v] (b0) at (2,1) {};
\node[v] (b1) at (0,0) {};
\node[v] (b2) at (1,1) {};
\node[v] (b3) at (1,0) {};
\node[v] (b4) at (1,-1) {};
\node[v] (b5) at (2,0) {};

\draw (b0)--(b2);
\draw (b0)--(b5);
\draw (b1)--(b3);
\draw (b1)--(b4);
\draw (b2)--(b3);
\draw (b3)--(b4);
\draw (b3)--(b5);
\end{scope}

\end{tikzpicture}
\caption{Laplacian cospectral, but not adjacency cospectral} 
\label{L-cospectral}
\end{center}
\end{minipage}
\end{center}
\end{figure}

There are 21 graphs on five vertices, and their graph spectra, both adjacent and Laplacian, are all distinct; consequently, both their $q$-spectra and transition $q$-spectra are also distinct. On the other hand, there are two triples and one pair of graphs with the same magnitude. 

There are 112 graphs on six vertices. There is exactly one pair of adjacent cospectral graphs (Figure \ref{cospectral-2}), but their $q$-spectra are distinct. There are two pairs of Laplacian cospectral graphs (Figure \ref{L-cospectral}), but their transition $q$-spectra are distinct. On the other hand, there are 20 classes comprising 63 graphs that share the same magnitude.

There are 853 graphs on seven vertices, and there are 11 pairs of graphs with identical $q$-spectra; all of these also have the same magnitude. In addition to these, there are 19 further pairs and one triple of adjacency cospectral graphs. 
There are 50 pairs and 5 triples of Laplacian cospectral graphs. All of them can be distinguished by the transition $q$-spectrum. 
As for the magnitude, there are 131 classes comprising 551 graphs that share the same magnitude. 

There are also 11 pairs having the same $d^{(1)}$-spectrum. Of these, 10 pairs are $q$-cospectral and one is not. 
We list the 10 pairs which are both $q$-cospectral and $d^{(1)}$-cospectral by their edge lists, and the remaining pairs by Figures \ref{q-cospectral11} and \ref{d1-cospectral}.

\begin{enumerate}
\item $\{12,13,14,15,24,35,46,47,56,57\}$ and $\{12,13,14,15,16,24,35,47,57,67\}$,
\item $\{12,13,14,15,16,24,35,45,46,47,56,57\}$ and $\{12,13,14,15,16,17,24,35,45,46,56,67\}$,
\item $\{12,13,14,15,16,24,26,35,47,57,67\}$ and $\{12,13,14,15,24,26,35,46,47,56,57\}$,
\item $\{12,13,14,15,17,23,26,34,45,56,67\}$ and $\{12,13,14,15,23,24,35,46,47,56,57\}$,
\item $\{12,13,14,15,16,23,24,35,36,47,57,67\}$ and $\{12,13,14,15,23,24,35,36,46,47,56,57\}$,
\item $\{12,13,14,15,16,23,24,35,36,37,45,56,57\}$ and $\{12,13,14,15,16,17,24,35,37,45,46,56,67\}$,
\item $\{12,13,14,15,24,26,35,36,46,47,56,57\}$ and $\{12,13,14,15,16,24,26,35,36,47,57,67\}$,
\item $\{12,13,14,15,16,24,27,35,37,45,46,56,57\}$ and $\{12,13,14,15,23,24,26,35,36,46,47,56,57\}$,
\item $\{12,13,14,15,16,17,23,24,26,35,36,47,57,67\}$ and $\{12,13,14,15,16,23,24,26,27,35,36,37,47,57\}$,
\item $\{12,13,14,15,16,17,23,24,26,35,37,46,47,56,57\}$ and $\{12,13,14,15,16,23,24,26,27,35,37,46,56,57,67\}$,
\item Figure \ref{q-cospectral11} ($q$-cospectral) and Figure \ref{d1-cospectral} ($d^{(1)}$-cospectral). 
\end{enumerate}

\begin{figure}[htbp]
\begin{center}
\begin{minipage}{.45\linewidth}
\begin{center}
\begin{tikzpicture}[
  scale=0.75,
  v/.style={circle, draw=black, fill=black, inner sep=1.4pt},
  e/.style={line width=0.8pt, draw=black}
]

% lower left graph
\begin{scope}[yshift=-4.0cm]
\node[v] (c1) at (0,0.0) {};
\node[v] (c2) at (1.05,0.6) {};
\node[v] (c3) at (1.05,-0.6) {};
\node[v] (c4) at (2.1,0) {};
\node[v] (c5) at (2.0,1.2) {};
\node[v] (c6) at (3.0,1.2) {};
\node[v] (c7) at (0,1.2) {};

\draw[e] (c1)--(c2)--(c4)--(c3);
\draw[e] (c1)--(c3)--(c2);
\draw[e] (c2)--(c5)--(c6);
\draw[e] (c2)--(c7);
\end{scope}

% lower right graph
\begin{scope}[xshift=5.6cm,yshift=-4.0cm]
\node[v] (d1) at (0,0.0) {};
\node[v] (d2) at (1.05,0.6) {};
\node[v] (d3) at (1.3,1.55) {};
\node[v] (d4) at (1.05,-0.6) {};
\node[v] (d5) at (2.1,0) {};
\node[v] (d6) at (2.9,0.8) {};
\node[v] (d7) at (2.9,-0.8) {};

\draw[e] (d1)--(d2)--(d3);
\draw[e] (d1)--(d4)--(d5)--(d2);
\draw[e] (d2)--(d4);
\draw[e] (d5)--(d6);
\draw[e] (d5)--(d7);
\end{scope}

\end{tikzpicture}
\caption{$q$-cospectral, not $d^{(1)}$-cospectral}
\label{q-cospectral11}
\end{center}
\end{minipage}
\hskip 0.8cm
\begin{minipage}{.45\linewidth}
\begin{center}
\begin{tikzpicture}[
  scale=0.7,
  v/.style={circle, draw=black, fill=black, inner sep=1.4pt},
  e/.style={line width=0.8pt, draw=black}
]

% upper left graph
\begin{scope}
\node[v] (a1) at (0,0.8) {};
\node[v] (a2) at (1.2,1.6) {};
\node[v] (a3) at (2.4,2.0) {};
\node[v] (a4) at (2.4,0.8) {};
\node[v] (a5) at (1.2,0.0) {};
\node[v] (a6) at (3.6,0.8) {};
\node[v] (a7) at (2.4,-0.4) {};

\draw[e] (a1)--(a2)--(a3)--(a6)--(a7)--(a5)--(a1);
\draw[e] (a2)--(a5);
\draw[e] (a2)--(a4);
\draw[e] (a3)--(a5);
\draw[e] (a4)--(a5);
\draw[e] (a4)--(a6);
\end{scope}

% upper right graph
\begin{scope}[xshift=5.6cm]
\node[v] (b1) at (0,0.8) {};
\node[v] (b2) at (0.9,1.7) {};
\node[v] (b3) at (2.4,1.7) {};
\node[v] (b4) at (3.4,0.85) {};
\node[v] (b5) at (2.4,0.0) {};
\node[v] (b6) at (1.4,0.567) {};
\node[v] (b7) at (0.9,0.0) {};

\draw[e] (b1)--(b6);
\draw[e] (b1)--(b7);
\draw[e] (b2)--(b3)--(b4)--(b5)--(b7)--(b2);
\draw[e] (b2)--(b6);
\draw[e] (b3)--(b5);
\draw[e] (b3)--(b6);
\draw[e] (b5)--(b6);
\end{scope}

\end{tikzpicture}
\caption{$d^{(1)}$-cospectral, not adjacency cospectral hence not $q$-cospectral}
\label{d1-cospectral}
\end{center}
\end{minipage}
\end{center}
\end{figure}

%!!!!!!!!!!!!!!!!!!!!!!!!!!!!!!!!!!!!!%%%%%%%%%%%%%%%%%%%%%%%%%%%%%%%%%%%%
\section{Open problems}\label{section_last}
%!!!!!!!!!!!!!!!!!!!!!!!!!!!!!!!!!!!!!%%%%%%%%%%%%%%%%%%%%%%%%%%%%%%%%%%%%

It remains unclear whether, in terms of distinguishing power, the transition $q$-spectrum is stronger than the $q$-spectrum, and both the $q$-spectrum and the transition $q$-spectrum are stronger than magnitude, or whether no strict ordering exists among them. 

%!!!!!!!!!!!!!!!!!!!!!!!!!!!!!!!!!!!!!%%%%%%%%%%%%%%%%%%%%%%%%%%%%%%%%%%%%
\section{Appendix}\label{section_appendix}
%!!!!!!!!!!!!!!!!!!!!!!!!!!!!!!!!!!!!!%%%%%%%%%%%%%%%%%%%%%%%%%%%%%%%%%%%%

We provide the details of Step 2 in the proof of Theorem {\rm \bf \ref{thm_four-point}}, postponed from Subsection \ref{subsection_four_pts}, that four-point metric spaces are determined by their $q$-spectra.

\begin{definition} \rm 
By {\em data of combination of opposite sides} (abbreviated as {\em dcos.}) for $\mathcal{S}_1=[a,b,c,d,f,g]$ and $\mathcal{S}_{\rm opp}=[\a,\beta,\ga]$ we mean a multiset $[\,[a,d],[b,f],[c,g]\,]$ that satisfies $[\a,\beta,\ga]=[a+d,b+f,c+g]$. 
We denote it by $\mathcal{C}$. 
\end{definition}

For any {\sl dcos} $\mathcal{C}$, there are at most two isometry classes that realize $\mathcal{C}$. 

\begin{figure}[htbp]
\begin{center}
\includegraphics[width=.45\linewidth]{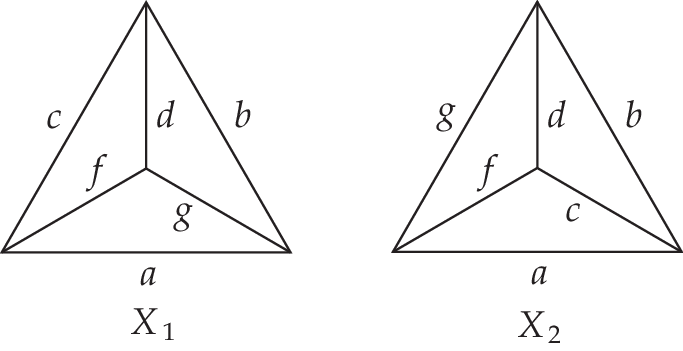}
\caption{}
\label{tetrahedra}
\end{center}
\end{figure}

(2-i) Suppose both $X$ and $X'$ with the same $\mathcal{S}_1$ and $\mathcal{S}_{\rm opp}$ have the same {\sl dcos}. 
We may assume without loss of generality that it is given by $[\,[a,d],[b,f],[c,g]\,]$. 
Then there are two possibilities as illustrated in Figure \ref{tetrahedra} to be fixed to identify $X$. 
Let $\varphi_i=\varphi_i(a,b,c,d,f,g)$ (or $\psi_j$) $(1\le i,j\le 4)$ be the perimeters of four triangles of $X_1$ (or resp. $X_2$); 
\[\begin{array}{ll}
\varphi_1(a,b,c,d,f,g)=a+b+c, &\quad \psi _1(a,b,c,d,f,g)=a+b+g, \\
\varphi_2(a,b,c,d,f,g)=a+f+g, &\quad \psi _2(a,b,c,d,f,g)=a+c+f, \\
\varphi_3(a,b,c,d,f,g)=b+d+g, &\quad \psi _3(a,b,c,d,f,g)=b+c+d, \\
\varphi_4(a,b,c,d,f,g)=c+d+f, &\quad \psi _4(a,b,c,d,f,g)=d+f+g. 
\end{array}\]
Suppose $\mathcal{S}_3^\triangle(X_1)=\mathcal{S}_3^\triangle(X_2)$. 
There are $24$ possibilities for the system of linear equations obtained from  $\mathcal{S}_3^\triangle(X_1)=\mathcal{S}_3^\triangle(X_2)$, each of which implies that $X_1$ is isometric to $X_2$. 
To be precise, for each permutation $\si\in\mathfrak{S}_4$ of four letters $\{1,2,3,4\}$, a system of linear equations 
\[
\left\{\begin{array}{l}
\varphi_1(a,b,c,d,f,g)=\psi_{\si(1)}(a,b,c,d,f,g), \\
\varphi_2(a,b,c,d,f,g)=\psi_{\si(2)}(a,b,c,d,f,g), \\
\varphi_3(a,b,c,d,f,g)=\psi_{\si(3)}(a,b,c,d,f,g), \\
\varphi_4(a,b,c,d,f,g)=\psi_{\si(4)}(a,b,c,d,f,g)
\end{array}\right.
\]
implies that at least one of $a=d$, $b=f$ or $c=g$ holds\footnote{The author used Maple 2023 for the computation of systems of linear equations in this subsection.}. 

(2-ii) Suppose $X$ and $X'$ with the same $\mathcal{S}_1$ and $\mathcal{S}_{\rm opp}$ have different {\sl dcos}. 
Then, at least one of the elements of $\mathcal{S}_{\rm opp}$, $\a,\beta,\ga$, say $\a$, is decomposed in two ways, say
\[\a=a+b=(a-\e)+(b+\e) \qquad (\e\ne0, \> a-\e\ne b).\]
Let $c$ and $d$ be the remaining elements of $\mathcal{S}_1$, i.e., 
$\mathcal{S}_1=[a,b,a-\e,b+\e,c,d].$ 
Put $\de=c+d$. 
Then $\de\ne\a$ since if $\de=\a$ then there is only one {\sl dcos}, which is a contradiction. 

\begin{lemma}\label{lem_docs_aab}
None of $X$ and $X'$ has {\sl dcos} of the form 
$\mathcal{C}_{\a\a\de}=[\,[a,b], [a-\e, b+\e], [c,d]\,].$ 
\end{lemma}

\begin{proof}
Assume one of $X$ and $X'$ has {\sl dcos} of the form $[\,[a,b], [a-\e, b+\e], [c,d]\,].$ 
Then $\mathcal{S}_{\rm opp}=[\a,\a,\de]$. 
Let $\mathcal{C}'$ be another {\sl dcos}; 
\[
\mathcal{C}'\colon
\a=A_{11}+A_{12}, \, \a=A_{21}+A_{22}, \, \de=B_1+B_2.
\]
Then at least two of $A_{ij}$ $(1\le i,j\le 2)$ belong to $\{a,b,a-\e,b+\e\}$. 
If the two are $A_{i1}$ and $A_{i2}$, as $[A_{(3-i)1},A_{(3-i)2}]$ cannot be $[c,d]$ because $\de\ne\a$, $\mathcal{C}'$ coincides with $\mathcal{C}_{\a\a\de}$, which is a contradiction. 
If the two are $A_{1j}$ and $A_{2j'}$ then $\mathcal{C}'$ coincides with $\mathcal{C}_{\a\a\de}$, which is a contradiction. 
\end{proof}

Now we may assume without loss of generality $X$ has {\sl dcos} 
\[\mathcal{C}_{\rm ii,0}\colon \a=(a-\e)+(b+\e), \, \beta=a+c, \,\ga=b+d.\]
Lemma \ref{lem_docs_aab}, $c+d\ne\a$, $\e\ne0$, and $a-\e\ne b$ imply that the {\sl dcos} of $X'$ is either 
\[\mathcal{C}_{\rm ii,1}\colon \a=a+b, \, \beta=(a-\e)+d, \,\ga=(b+\e)+c\]
or 
\[\mathcal{C}_{\rm ii,2}\colon \a=a+b, \, \beta=(b+\e)+d, \,\ga=(a-\e)+c.\]

\begin{figure}[htbp]
\begin{center}
\begin{minipage}{.45\linewidth}
\begin{center}
\includegraphics[width=\linewidth]{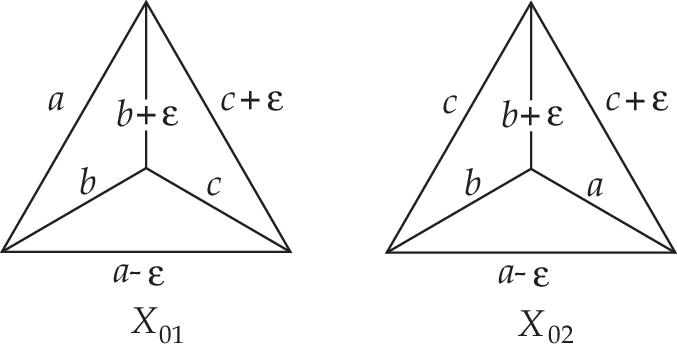}
\caption{}
\label{tetrahedra_0102}
\end{center}
\end{minipage}
\hskip 0.8cm
\begin{minipage}{.45\linewidth}
\begin{center}
\includegraphics[width=\linewidth]{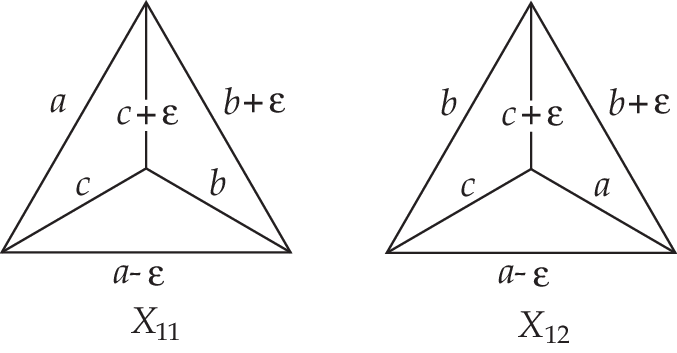}
\caption{}
\label{tetrahedra_1112}
\end{center}
\end{minipage}
\end{center}
\end{figure}

(2-ii-1) Suppose $X'$ has {\sl dcos} $\mathcal{C}_{\rm ii,1}$. 
Then $\beta=a+c=a+d-\e$, which implies that $d=c+\e$, and hence $\mathcal{C}_{\rm ii,0}$ and $\mathcal{C}_{\rm ii,1}$ are given by 
\[
\begin{array}{lll}
\mathcal{C}_{\rm ii,0}\colon \a=(a-\e)+(b+\e), &\, \beta=a+c, &\,\ga=b+(c+\e), \\
\mathcal{C}_{\rm ii,1}\colon \a=a+b, &\, \beta=(a-\e)+(c+\e), &\,\ga=(b+\e)+c.
\end{array}
\]
There are two realizations of $\mathcal{C}_{\rm ii,0}$, $X_{01}$ and $X_{02}$ (Figure \ref{tetrahedra_0102}), and two realizations of $\mathcal{C}_{\rm ii,1}$, $X_{11}$ and $X_{12}$ (Figure \ref{tetrahedra_1112}). 
Suppose $X_{0i}$ and $X_{1j}$ $(i,j=1,2)$ have the same $\mathcal{S}_3^\triangle$. 
Then at least one of the systems of four linear equations must be satisfied. 
Examine all the $4\times 4!=96$ systems of linear equations, and we obtain either the isometry between $X_{0i}$ and $X_{1j}$ or $\e=0$, which contradicts our assumption (2-ii). 
Let us give two examples. 
Assume $X_{01}$ and $X_{11}$ have the same $\mathcal{S}_3^\triangle$. 
Put 
\[\begin{array}{ll}
\varphi_1(a,b,c,\e)=2a+c, &\quad \psi _1(a,b,c,\e)=2a+b, \\
\varphi_2(a,b,c,\e)=a+b+c-\e, &\quad \psi _2(a,b,c,\e)=a+b+c-\e, \\
\varphi_3(a,b,c,\e)=a+2b+\e, &\quad \psi _3(a,b,c,\e)=a+2c+\e, \\
\varphi_4(a,b,c,\e)=b+2c+2\e, &\quad \psi _4(a,b,c,\e)=2b+c+2\e. 
\end{array}\]

\begin{figure}[htbp]
\begin{center}
\includegraphics[width=.45\linewidth]{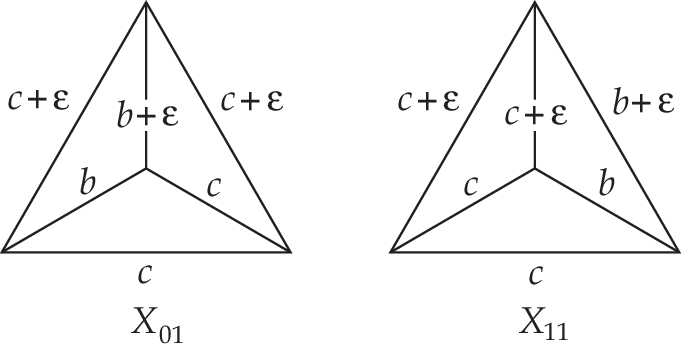}
\caption{}
\label{ii-example1}
\end{center}
\end{figure}

$\bullet$ The system of linear equations $\varphi_1=\psi_3, \varphi_2=\psi_2, \varphi_3=\psi_4, \varphi_4=\psi_1$, i.e. 
\[
\left\{\begin{array}{rl}
2a+c &=a+2c+\e,  \\
a+b+c-\e &=a+b+c-\e,  \\
a+2b+\e &=2b+c+2\e,  \\
b+2c+2\e &=2a+b 
\end{array}\right.
\]
has a solution $a=c+\e$. Then $X_{01}$ and $X_{11}$ are as illustrated in Figure \ref{ii-example1}, and they are isometric. 

$\bullet$ The system of linear equations $\varphi_1=\psi_1, \varphi_2=\psi_3, \varphi_3=\psi_2, \varphi_4=\psi_4$, i.e. 
\[
\left\{\begin{array}{rl}
2a+c &=2a+b,  \\
a+b+c-\e &=a+2c+\e,  \\
a+2b+\e &=a+b+c-\e,  \\
b+2c+2\e &=2b+c+2\e 
\end{array}\right.
\]
has a solution $b=c$ and $\e=0$, which contradicts our assumption that $\e\ne0$. 

\smallskip
(2-ii-2) The case when $X'$ has {\sl dcos} $\mathcal{C}_{\rm ii,2}$ can similarly be proved by exhaustion.

\noindent
Jun O'Hara

\noindent
Department of Mathematics and Informatics, Faculty of Science, 
Chiba University

\noindent
1-33 Yayoi-cho, Inage, Chiba, 263-8522, JAPAN.  

\noindent
E-mail: ohara@math.s.chiba-u.ac.jp

\end{document}